\documentclass{amsart}

\usepackage{latexsym}
\usepackage{amsmath}
\usepackage{amssymb}

\begin{document}



\newtheorem{theorem}{Theorem}[section]
\newtheorem{corollary}[theorem]{Corollary}
\newtheorem{proposition}[theorem]{Proposition}
\newtheorem{lemma}[theorem]{Lemma}
\newtheorem{conjecture}[theorem]{Conjecture}

\theoremstyle{remark}
\newtheorem{remark}[theorem]{Remark}

\theoremstyle{definition}
\newtheorem{definition}[theorem]{Definition}
\newtheorem{example}[theorem]{Example}

\numberwithin{equation}{section}
\numberwithin{theorem}{section}


\newcommand{\Iz}{I_{Z}}
\newcommand{\Ix}{I_X}
\newcommand{\CIgr}{CI_{grid}}
\newcommand{\CIgen}{CI_{gen}}
\newcommand{\Xgr}{X_{grid}}
\newcommand{\Xgen}{X_{gen}}
\newcommand{\Ygr}{Y_{grid}}
\newcommand{\Ygen}{Y_{gen}}
\newcommand{\Zgr}{Z_{grid}}
\newcommand{\Zgen}{Z_{gen}}
\newcommand{\Iy}{I_Y}
\newcommand{\A}{\mathcal{A}}
\newcommand{\C}{\mathcal{C}}
\newcommand{\U}{\mathcal{U}}
\newcommand{\Y}{\mathbb{Y}}
\newcommand{\s}{\mathbb{S}}
\newcommand{\Z}{\mathbb{Z}}
\newcommand{\Zr}{\mathbb{Z}_{red}}
\newcommand{\N}{\mathbb{N}}
\newcommand{\pr}{\mathbb{P}}
\newcommand{\X}{\mathbb{X}}
\newcommand{\supp}{\operatorname{Supp}}
\newcommand{\HH}{\mathcal{H}}
\newcommand{\M}{\mathcal{M}}
\newcommand{\G}{\mathcal{G}}
\newcommand{\F}{\mathcal{F}}
\newcommand{\K}{\mathcal{K}}
\newcommand{\LL}{\mathcal{L}}
\newcommand{\rk}{\operatorname{rk}}
\newcommand{\ri}{\operatorname{ri}}
\newcommand{\ox}{\overline{x}}



\title[Powers of complete intersections]{Powers
of complete intersections: Graded Betti numbers and
applications}
\thanks{Revised Version: March 17, 2005}

\author{Elena Guardo}
\address{Dipartimento di Matematica e Informatica\\
Viale A. Doria, 6 - 95100 - Catania, Italy}
\email{guardo@dmi.unict.it}
\author{Adam Van Tuyl}
\address{Department of Mathematical Sciences \\
Lakehead University \\
Thunder Bay, ON P7B 5E1, Canada}
\email{avantuyl@sleet.lakeheadu.ca}




\begin{abstract}
Let $I = (F_1,\ldots,F_r)$ be a homogeneous ideal of the ring
$R = k[x_0,\ldots,x_n]$ generated by a regular sequence of
type $(d_1,\ldots,d_r)$.  
We give an elementary proof for an explicit description of the 
graded Betti numbers of $I^s$  for any $s \geq 1$.
These numbers depend only upon the type and $s$.
We then use this description
to: (1) write $H_{R/I^s}$, the Hilbert function of $R/I^s$,
in terms of $H_{R/I}$; (2) verify that the $k$-algebra
$R/I^s$ satisfies a conjecture of Herzog-Huneke-Srinivasan; and
(3) obtain information about the numerical invariants 
associated to sets of fat points in $\mathbb{P}^n$ whose support
is a complete intersection or a complete intersection
minus a point.
\end{abstract}
\maketitle


\section*{Introduction}
In this paper we give an explicit description
of the graded
Betti numbers of a power of a complete intersection and provide
some applications of this result.

It is well known that the graded minimal free resolution of a homogeneous
complete intersection $I = (F_1,\ldots,F_r) \subseteq R = k[x_0,\ldots,x_n]$
is given by the Koszul resolution.    
Furthermore, the graded Betti numbers of $I$ depend only upon the type $(d_1,\ldots,d_r)$
where $d_i = \deg F_i$ (cf. Theorem \ref{ciproperties}).  
However, the description of the graded Betti numbers of a power of a 
complete intersection does not enjoy the same level of familiarity. 

It has long been known that $I^s$, the power
of a complete intersection, can be realized as the ideal
generated by the maximal minors of a $s \times (s+r-1)$ matrix with
entries consisting of the $F_i$'s.  Since $I^s$ is an example
of a determinantal ideal, a minimal free resolution of $I^s$
can be obtained by using the Eagon-Northcott complex.  However, the
maps in this resolution may not be degree preserving, so the graded
Betti numbers cannot be extracted.

Alternatively, a resolution of $I^s$ can be found in \cite{BE}
(see \cite[Theorem 2.1]{S} for a more relevant formulation).
Implicit in \cite{S} is the fact that 
the minimal resolution of $I^s$ is a graded minimal resolution of $I^s$
provided $I$ is homogeneous.  The graded Betti
numbers of $I^s$ are a consequence of this result, but
to the best of our knowledge, they have never been explicitly
described.  This probably accounts for the lack of familiarity mentioned
above for these invariants.

Our main result (Theorem \ref{maintheorem}) is to give an explicit
description for the graded Betti numbers of $I^s$. As we shall show,
these numbers can be computed directly from
 the type $(d_1,\ldots,d_r)$ and $s$.
To further differentiate our result from past work, we
give an elementary proof which avoids the machinery used
in \cite{BE,S}.  
We relate the ideals $I^s$, $I^{s-1}$ and
$(F_2,\ldots,F_r)^s$ through a basic double link.  Then, by using a mapping cone construction
and induction on the tuple $(r,s)$, we obtain the graded Betti
numbers of $I^s$.

As a corollary, we express $H_{R/I^s}$, the Hilbert function
of $R/I^s$, as a function of $H_{R/I}$.  Like the graded Betti
numbers of $I^s$, we were unable to find in the literature an expression
for $H_{R/I^s}$.
This omission is  also noted in \cite[Page 799]{GGH}.   


The final three sections are devoted to applications of Theorem
\ref{maintheorem}.  In \S 3, we verify that the
$k$-algebra $R/I^s$ satisfies a conjecture of Herzog-Huneke-Srinivasan which
relates the multiplicity $e(R/I^s)$ to the shifts in the graded minimal
free resolution of $R/I^s$.  

In \S 4, we use Theorem \ref{maintheorem}
to study sets of fat points $Z \subseteq \pr^n$ whose support is a complete
intersection.  When $Z$ is homogeneous (all the points have the same
multiplicity), then Theorem \ref{maintheorem} gives the graded
minimal free resolution of $I_Z$.  When $Z$ is not homogeneous, 
we obtain partial information about the invariants associated
to $I_Z$.  

In \S 5, we investigate sets of 
fat points in $\pr^n$ whose support is a complete
intersection minus a point.  We give examples to show
different constructions of the underlying complete 
intersection, e.g. whether or not the defining forms are irreducible
or reducible, may result in different numerical invariants, thus 
 implying extra hypotheses on
the support are needed.  
We therefore study such schemes under
the extra condition that the underlying
complete intersection can ``split'' into smaller complete 
intersections.  This restriction allows us to use
Theorem \ref{maintheorem} to obtain bounds
on the least degree of a form passing through the scheme and
the regularity index.  This section extends
some of the results of \cite{BGVT1,BGVT2,Gu1} which studied this question
in $\pr^2$.

\noindent
{\it Acknowledgments.} We wish to thank
B. Harbourne, C. Peterson, H. Srinivasan, and J. Weyman for comments and for
answering some of our questions during the early stages of this
paper.  We especially thank J. Migliore for  bringing to our
attention \cite{KMMNP}, which lead to Lemma \ref{inductionstep},
thereby greatly simplifying our original arguments.  The
program {\tt CoCoA} \cite{C} was used extensively during the
preliminary stages. 
The second
author also acknowledges the financial support of NSERC.


\section{Complete intersections}

Let $k$ denote an infinite
field with char$(k) = 0$, and set $R = k[x_0,\ldots,x_n]$.
A homogeneous ideal $I = (F_1,\ldots,F_r) \subseteq R$ is a
{\it complete intersection} if $F_1,\ldots, F_r$ form a regular
sequence on $R$.  If $\deg F_i = d_i$, then we say that $I$ is a complete
intersection of type $(d_1,\ldots,d_r)$.  
 Since $F_1,\ldots,F_r$ are homogeneous,
any permutation of their order again results in a regular sequence.
Thus, we can assume that $d_1 \leq d_2 \leq \cdots \leq d_r$.
The graded Betti numbers in the graded minimal free resolution of $I$ 
depend only upon the type as described below.

\begin{theorem}[Koszul Resolution] \label{ciproperties}
Let $I \subseteq R$ be a complete intersection
of type $(d_1,\ldots,d_r)$.
Then the graded minimal free resolution of $I$ has the form
\[0 \rightarrow \F_{r-1} \rightarrow \F_{r-2}
\rightarrow \cdots \rightarrow \F_1 \rightarrow \F_0 \rightarrow I \rightarrow
0\] where
\[{\displaystyle
\F_j = \bigoplus_{1 \leq i_1 < i_2
< \cdots < i_{j+1} \leq r} R(-d_{i_1}-d_{i_2}- \cdots- d_{i_{j+1}})}
~~\mbox{for  $j = 0,\ldots,r-1$.}
\]
\end{theorem}

Some of the properties of the ideal $I^s$, where $I$ is a complete
intersection and $s \in \N^+$, are given below.

\begin{theorem} \label{determinate}
Let $I = (F_1,\ldots,F_r)$
be an ideal of $R = k[x_0,\ldots,x_n]$ generated by a regular sequence
of length $r \leq  n+1$, and let $\mathcal{A}$ be
the $s \times (s+r-1)$ matrix
\[ \mathcal{A} = \begin{pmatrix}
F_1 & F_2 & \cdots & F_r & 0 & 0 &0  &\cdots & 0 \\
0 & F_1 & F_2 & \cdots & F_r & 0 &0  & \cdots & 0 \\
0 & 0 & F_1 & F_2 & \cdots & F_r & 0 & \cdots & 0 \\
\cdot & \cdot&\cdot & \cdot&\cdot & \cdot&\cdot & \cdots& \cdot \\
0 & 0 & 0& \cdots & 0& F_1 & F_2 & \cdots & F_r \\
\end{pmatrix}.\]
Then
\begin{enumerate}
\item[$(i)$] $I^s = I_s(\mathcal{A})$, the ideal generated by
the $s \times s$ maximal minors of $\mathcal{A}$.
\item[$(ii)$] for each $s \geq 1$, the ideal $I^s$ is perfect
of grade $r$.
\item[$(iii)$] the minimal free
resolution of $I^s$ is given by the Eagon-Northcott
complex constructed from $\mathcal{A}$.
\end{enumerate}
\end{theorem}

\begin{proof}
Statement $(i)$ is in \cite{BV} following Remark 2.13.
Statement $(ii)$ is \cite[Proposition 2.14]{BV}.  
Since the grade of $I_s(\mathcal{A})$ is  $r$,
the Eagon-Northcott complex constructed from $\mathcal{A}$
gives a minimal free resolution of $R/I_s(\mathcal{A})$
by \cite[Theorem 5.2]{EN}.
\end{proof}

Although the Eagon-Northcott complex gives a minimal resolution
of $I^s$, the maps in the resolution may not be degree
preserving, and thus, we cannot use this resolution
to compute the graded Betti numbers except under
extra hypotheses on the type, e.g., \cite{BN}.  To see this,
note that the matrix in the above theorem defines a map
$R^{s+r-1} \stackrel{\mathcal{A}}{\longrightarrow} R^{s}$.
This map may fail to be a homogeneous map, e.g.,
if the degrees of the $F_i$'s fail to be an arithmetic progression,
and thus the maps in Eagon-Northcott complex are not homogeneous.
However, we can use the Eagon-Northcott resolution of $I^s$
to determine the ranks of the modules in the
 resolution, thus allowing us to count the number
of generators of each syzygy module.

\begin{corollary}   \label{ranks}
With the hypotheses as in the previous theorem, let
\[ 0 \rightarrow \mathcal{F}_{r-1} \rightarrow \mathcal{F}_{r-2}
\rightarrow \cdots \rightarrow \mathcal{F}_0 \rightarrow I^s \rightarrow
0\]
be the minimal free resolution of $I^s$ given by the
Eagon-Northcott complex.  Then
\[\operatorname{rk} \mathcal{F}_i = \binom{r+s-1}{s+i}\binom{s-1+i}{i} ~~
\mbox{for $i = 0,\ldots,r-1$}.\]
\end{corollary}

The following lemma relates $I^s$ to a complete intersection
with one less generator.

\begin{lemma}   \label{inductionstep}
Let $F_1,\ldots,F_r$ be a regular sequence in $R$
with $\deg F_i = d_i$.  Set $I = (F_1,\ldots,F_r)$ and
$J = (F_2,\ldots,F_r)$.  Then for each positive integer $s$,
\[I^s = J^{s} + F_1\cdot I^{s-1}.\]
Furthermore, we have the following short exact sequence
\begin{equation} \label{excseq}
0 \longrightarrow J^s(-d_1){\longrightarrow} I^{s-1}(-d_1)\oplus
J^{s}{\longrightarrow} J^{s}+F_1\cdot I^{s-1} = I^s\longrightarrow 0.
\end{equation}
\end{lemma}

\begin{proof}
The equality of ideals is immediate.
The short exact sequence is from \cite[Lemma 4.8]{KMMNP}.
\end{proof}


\section{Graded Betti numbers of powers of complete intersections}

Let $I =(F_1,\ldots,F_r)$ be a homogeneous complete intersection of
 $R = k[x_0,\ldots,x_n]$.
For each $s \in \N^+$ we describe how the graded
Betti numbers in the graded minimal free resolution of $I^s$
depend only upon the type $(d_1,\ldots,d_r)$ and $s$.
To describe the resolution we introduce
the sets
\[\mathcal{M}_{r,s,t} := \left\{ (a_1,\ldots,a_r)\in \N^r
\left|\begin{array}{l}
a_1 + \cdots + a_r = s ~~\mbox{and at}\\
\mbox{least $t$ of the $a_i$'s are non-zero}
\end{array} \right\}\right.\]
for all positive integers $r,s,$ and $t$.
With this notation we have

\begin{theorem}     \label{maintheorem}
Let $I \subseteq R$ be a complete intersection 
of type $(d_1,\ldots,d_r)$ and let $s \in \N^+$.
Then the graded minimal free resolution of $I^s$ has the form
\[
0 \rightarrow \HH_{r-1} \rightarrow \HH_{r-2} \rightarrow
\cdots \rightarrow \HH_{0} \rightarrow I^s \rightarrow 0
\] where
\[
\HH_0  =
\bigoplus_{(a_1,\ldots,a_r) \in \M_{r,s,1}} R(-a_1d_1-\cdots-a_rd_r)
\]
and for $i=1,\ldots,r-1$,
\footnotesize
\[
\HH_i  =  \bigoplus_{l_1 = i+1}^r
\left[\bigoplus_{l_2=l_1}^{r}\left[\cdots
\left[ \bigoplus_{l_i=l_{i-1}}^r
\left[\bigoplus_{(a_1,\ldots,a_r) \in \M_{r,s+i,l_i}}
R(-a_1d_1-\cdots-a_rd_r)\right]\right]\cdots\right]\right].
\]
\normalsize
\end{theorem}

\begin{proof}
The proof is by induction on the tuple $(r,s)$.  If $s =1$,
and $r \leq n+1$ is any positive integer, then the resolution of $I$
is given by the Koszul resolution (Theorem \ref{ciproperties}).  
The reader can verify that
the sets $\M_{r,s,t}$ account for the expected Betti numbers in this
case.  If $r=1$, and $s$ is any integer, then $I^s$ is
principal, and the result also follows.

So, let $(r,s)$ be a tuple with $1 < r \leq n+1$ and $s\in\N^+$, 
and assume that the result
holds for all ideals of the form $(F_1,\ldots,F_{r'})^{s'}$ with
$(r',s') <_{lex} (r,s)$ with respect to the lexicographical
ordering.  Set $I = (F_1,\ldots,F_r)$ and $J = (F_2,\ldots,F_r)$.
The short
exact sequence (\ref{excseq}) of Lemma \ref{inductionstep} 
relates $I^s, I^{s-1}$ and $J^s$:
\[
0 \longrightarrow J^s(-d_1){\longrightarrow} I^{s-1}(-d_1)\oplus
J^{s}{\longrightarrow} I^s
\longrightarrow 0.
\]
Let
\[
0 \rightarrow \G_{r-2} \rightarrow \G_{r-3} \rightarrow
\cdots \rightarrow \G_{0} \rightarrow J^s\rightarrow 0 ~~~~~~~~\text{and}
\]
\[
0 \rightarrow \K_{r-1} \rightarrow \K_{r-2} \rightarrow
\cdots \rightarrow \K_{0} \rightarrow I^{s-1}\rightarrow 0
\]
be the minimal free resolutions of $J^s$ and  $I^{s-1}$,
respectively.  The mapping
cone construction and the exact sequence
then gives the following resolution of $I^s$:
\begin{equation} \label{resolution}
0 \rightarrow \HH_{r-1} \rightarrow \HH_{r-2} \rightarrow
\cdots \rightarrow \HH_{0} \rightarrow I^s\rightarrow 0
\end{equation}
where
\begin{eqnarray*}
\HH_0 &= & \G_0\oplus \K_0(-d_1) \\
\HH_i &= & \G_{i-1}(-d_1)\oplus \G_i\oplus \K_i(-d_1) ~\mbox{for
$i=1,\cdots,r-2$} \\
\HH_{r-1} & = & \G_{r-2}(-d_1)\oplus \K_{r-1}(-d_1).
\end{eqnarray*}

We verify first that $(\ref{resolution})$ is minimal.  The resolution
has the correct length by Theorem \ref{determinate}.  To show that
the sequence does not split, it is enough to apply induction
to show that the rank
of each $\HH_i$ equals the rank expected by Corollary ~\ref{ranks}.
Indeed, for $i = 1,\ldots,r-2$ we get
\footnotesize
\begin{eqnarray*}
\rk \HH_i & = & \rk \G_{i-1} + \rk \G_i + \rk \K_i\\
& = & \binom{r+s-2}{s+i-1}\binom{s+i-2}{i-1} +
\binom{r+s-2}{s+i}\binom{s+i-1}{i}+
\binom{r+s-2}{s+i-1}\binom{s+i-2}{i}\\
& = &\binom{r+s-1}{s+i}\binom{s-i-1}{i}.
\end{eqnarray*}
\normalsize
The proofs that $\HH_0$ and $\HH_{r-1}$ have the correct ranks follow
similarly, so we have omitted them.  It  now follows that $(\ref{resolution})$
is minimal.

We now compute the graded Betti numbers of $I^s$ by applying the
induction hypothesis to $J^s$ and $I^{s-1}$.
Due to the tedious nature of the proof, we shall only 
show that $\HH_{i}$ has the correct
Betti numbers for $i = 2,\ldots,r-2$, and omit
the similar proofs for $\HH_0,\HH_1,$ and $\HH_r$.

For $i=2,\cdots,r-2$ we have
$\HH_i  =  \G_{i-1}(-d_1)\oplus \G_i\oplus \K_i(-d_1)$.  The induction
hypothesis applied to $\G_{i-1}(-d_1)$ gives
\allowdisplaybreaks
\footnotesize
\begin{eqnarray*}
\G_{i-1}(-d_1) & = & \bigoplus_{l_2=i}^{r-1}
\left[\bigoplus_{l_3=l_2}^{r-1}\left[
\cdots\left[
\bigoplus_{l_i=l_{i-1}}^{r-1}
\left[
\bigoplus_{\M_{r-1,s+i-1,l_i}}R(-a_2d_2-\cdots-a_rd_r-d_1)
\right]\right]\cdots\right] \right]\\
& = & \bigoplus_{l_2=i}^{r-1}
\left[\bigoplus_{l_3=l_2}^{r-1}\left[
\cdots\left[
\bigoplus_{l_i=l_{i-1}}^{r-1}
\left[
\bigoplus_{\footnotesize{\begin{array}{c}
\M_{r,s+i,l_i+1}\\
\text{and}~a_1 = 1\end{array}}}
R(-a_1d_1 -\cdots-a_rd_r)
\right]\right]\cdots\right] \right]\\
& = & \bigoplus_{l_2=i+1}^{r}
\left[\bigoplus_{l_3=l_2}^{r}\left[
\cdots\left[
\bigoplus_{l_i=l_{i-1}}^{r}
\left[
\bigoplus_{\footnotesize{\begin{array}{c}
\M_{r,s+i,l_i}\\
\text{and}~a_1 = 1\end{array}}}
R(-a_1d_1 -\cdots-a_rd_r)
\right]\right]\cdots\right] \right].
\end{eqnarray*}
\normalsize

We can write $\G_i$ as
\footnotesize
\begin{eqnarray*}
\G_i & = &
 \bigoplus_{l_1=i+1}^{r-1}
\left[\bigoplus_{l_2=l_1}^{r-1}\left[
\cdots\left[
\bigoplus_{l_i=l_{i-1}}^{r-1}
\left[
\bigoplus_{\M_{r-1,s+i,l_i}}R(-a_2d_2-\cdots-a_rd_r)
\right]\right]\cdots\right] \right]\\
& = &  \bigoplus_{l_2=i+1}^{r-1}
\left[\bigoplus_{l_3=l_2}^{r-1}\left[
\cdots\left[
\bigoplus_{l_i=l_{i-1}}^{r-1}
\left[
\bigoplus_{\footnotesize{\begin{array}{c}
\M_{r,s+i,l_i}\\
\text{and}~a_1 = 0\end{array}}}R(-a_1d_1-\cdots-a_rd_r)
\right]\right]\cdots\right] \right] \oplus\\
&   & \hspace{.5cm}\vdots \\
&   &  \bigoplus_{l_2=r-1}^{r-1}
\left[\bigoplus_{l_3=l_2}^{r-1}\left[
\cdots\left[
\bigoplus_{l_i=l_{i-1}}^{r-1}
\left[
\bigoplus_{\footnotesize{\begin{array}{c}
\M_{r,s+i,l_i}\\
\text{and}~a_1 = 0\end{array}}}R(-a_1d_1-\cdots-a_rd_r)
\right]\right]\cdots\right] \right].\\
\end{eqnarray*}
\normalsize
Because the set $\{(a_1,\ldots,a_r)\in\M_{r,s+i,r} ~|~ a_1 = 0\} = \emptyset$,
$\G_i$ is equal to
\footnotesize
\begin{eqnarray*}
\G_i & = & \bigoplus_{l_2=i+1}^{r}
\left[\bigoplus_{l_3=l_2}^{r}\left[
\cdots\left[\bigoplus_{l_i=l_{i-1}}^{r}
\left[ \bigoplus_{\footnotesize{\begin{array}{c}\M_{r,s+i,l_i}\\
\text{and}~a_1 = 0\end{array}}}R(-a_1d_1-\cdots-a_rd_r)
\right]\right]\cdots\right] \right] \oplus\\
&   & \hspace{.5cm}\vdots \\
&   &  \bigoplus_{l_2=r-1}^{r}
\left[\bigoplus_{l_3=l_2}^{r}\left[
\cdots\left[
\bigoplus_{l_i=l_{i-1}}^{r}
\left[
\bigoplus_{\footnotesize{\begin{array}{c}
\M_{r,s+i,l_i}\\
\text{and}~a_1 = 0\end{array}}}R(-a_1d_1-\cdots-a_rd_r)
\right]\right]\cdots\right] \right] \oplus\\
&& \bigoplus_{l_2=r}^{r}
\left[\bigoplus_{l_3=l_2}^{r}\left[
\cdots\left[\bigoplus_{l_i=l_{i-1}}^{r}
\left[
\bigoplus_{\footnotesize{\begin{array}{c}
\M_{r,s+i,l_i}\\
\text{and}~a_1 = 0\end{array}}}R(-a_1d_1-\cdots-a_rd_r)
\right]\right]\cdots\right] \right].
\end{eqnarray*}
\normalsize
For $\K_i(-d_i)$, the induction hypothesis gives
\footnotesize
\begin{eqnarray*}
\K_i(-d_1) &=&
  \bigoplus_{l_1=i+1}^{r}
\left[\bigoplus_{l_2=l_1}^{r}\left[
\cdots\left[
\bigoplus_{l_i=l_{i-1}}^{r}
\left[
\bigoplus_{\M_{r,s+i-1,l_i}}R(-a_1d_1-\cdots-a_rd_r-d_1)
\right]\right]\cdots\right] \right]\\
& = &  \bigoplus_{l_2=i+1}^{r}
\left[\bigoplus_{l_3=l_2}^{r}\left[
\cdots\left[
\bigoplus_{l_i=l_{i-1}}^{r}
\left[
\bigoplus_{\M_{r,s+i-1,l_i}}R(-a_1d_1-\cdots-a_rd_r-d_1)
\right]\right]\cdots\right] \right] \oplus\\
& & \hspace{.5cm}\vdots \\
& &  \bigoplus_{l_2=r}^{r}
\left[\bigoplus_{l_3=l_2}^{r}\left[
\cdots\left[
\bigoplus_{l_i=l_{i-1}}^{r}
\left[
\bigoplus_{\M_{r,s+i-1,l_i}}R(-a_1d_1-\cdots-a_rd_r-d_1)
\right]\right]\cdots\right] \right].\\
\end{eqnarray*}
\normalsize
We have the following identity:
\begin{eqnarray*}
\bigoplus_{\M_{r,s+i-1,l_i}}R(-a_1d_1-\cdots-a_rd_r-d_1)
&=& \bigoplus_{\footnotesize{\begin{array}{c}
\M_{r,s+i,l_i}\\
\text{and}~a_1 \geq 2\end{array}}} R(-a_1d_1-\cdots-a_rd_r)
\oplus \\
&&
\bigoplus_{\footnotesize{\begin{array}{c}
\M_{r,s+i,l_i+1}\\
\text{and}~a_1 = 1 \end{array}}} R(-a_1d_1-\cdots-a_rd_r).
\end{eqnarray*}
Because of this identity, we can write $\K_i(-d_1)$ as
\footnotesize
\begin{eqnarray*}
\K_i(-d_1) &=& \bigoplus_{l_2=i+1}^{r}
\left[\bigoplus_{l_3=l_2}^{r}\left[
\cdots\left[
\bigoplus_{l_i=l_{i-1}}^{r}
\left[
\bigoplus_{\footnotesize{\begin{array}{c}
\M_{r,s+i,l_i}\\
\text{and}~a_1 \geq 2\end{array}}}R(-a_1d_1-\cdots-a_rd_r)
\right]\right]\cdots\right] \right] \oplus \\
&&\bigoplus_{l_2=i+2}^{r}
\left[\bigoplus_{l_3=l_2}^{r}\left[
\cdots\left[
\bigoplus_{l_i=l_{i-1}}^{r}
\left[
\bigoplus_{\footnotesize{\begin{array}{c}
\M_{r,s+i,l_i}\\
\text{and}~a_1 =1\end{array}}}R(-a_1d_1-\cdots-a_rd_r)
\right]\right]\cdots\right] \right] \oplus
\\
& &  \bigoplus_{l_2=i+2}^{r}
\left[\bigoplus_{l_3=l_2}^{r}\left[
\cdots\left[
\bigoplus_{l_i=l_{i-1}}^{r}
\left[
\bigoplus_{\footnotesize{\begin{array}{c}
\M_{r,s+i,l_i}\\
\text{and}~a_1 \geq 2\end{array}}}R(-a_1d_1-\cdots-a_rd_r)
\right]\right]\cdots\right] \right] \oplus\\
& & \hspace{.5cm}\vdots \\
& & \bigoplus_{l_2=r}^{r}
\left[\bigoplus_{l_3=l_2}^{r}\left[
\cdots\left[
\bigoplus_{l_i=l_{i-1}}^{r}
\left[
\bigoplus_{\footnotesize{\begin{array}{c}
\M_{r,s+i,l_i}\\
\text{and}~a_1 =1\end{array}}}R(-a_1d_1-\cdots-a_rd_r)
\right]\right]\cdots\right] \right] \oplus
\\
& &  \bigoplus_{l_2=r}^{r}
\left[\bigoplus_{l_3=l_2}^{r}\left[
\cdots\left[
\bigoplus_{l_i=l_{i-1}}^{r}
\left[
\bigoplus_{\footnotesize{\begin{array}{c}
\M_{r,s+i,l_i}\\
\text{and}~a_1 \geq 2\end{array}}}R(-a_1d_1-\cdots-a_rd_r)
\right]\right]\cdots\right] \right].
\end{eqnarray*}
\normalsize
As a consequence, since $\HH_i = \G_{i-1}(-d_1)\oplus
\G_i\oplus\K_i(-d_1)$, we have
\footnotesize
\[\HH_{i} = \bigoplus_{l_1 = i+1}^r
\left[\bigoplus_{l_2=l_1}^{r}\left[\cdots
\left[ \bigoplus_{l_i=l_{i-1}}^r
\left[\bigoplus_{(a_1,\ldots,a_r) \in \M_{r,s+i,l_i}}
R(-a_1d_1-\cdots-a_rd_r)\right]\right]\cdots\right]\right].
\]
\normalsize
as desired.
\end{proof}

\begin{remark}
The minimal resolution of $I^s$ as given in \cite[Theorem 2.1]{S}
is also graded if $I$ is homogeneous, although it is not explicitly 
stated as such.  For the 
convenience of the reader, we sketch out how to make
the resolution graded.  
Let $F$ be a free $R$-module with $\operatorname{rk} F = r$.  Let 
$L_i^sF$ denote the image of the map
$d_{i,s}: \bigwedge^i F \otimes S_{s-1}F 
\longrightarrow \bigwedge^{i-1} F \otimes
S_sF$
where $d_{i,s}$ is defined as in \cite[Definition 1.1]{S} and
$SF$ is the symmetric algebra on $F$ .  The
resolution of $I^s = (F_1,\ldots,F_r)^s$ then has the form:
\[0 \rightarrow L_r^sF \rightarrow L_{r-1}^sF \rightarrow \cdots \rightarrow
L_2^sF \rightarrow S_sF \rightarrow I^s \rightarrow 0.\]
As explained in \cite[Page 152]{S}, a
 basis for $ L_i^sF$ for each $i$ 
can be identified with the standard tableau of the form
\begin{center}
\begin{picture}(100,100)(0,0)
\put(0,0){\line(0,1){100}}
\put(0,0){\line(1,0){20}}
\put(20,0){\line(0,1){100}}
\put(40,80){\line(0,1){20}}
\put(80,80){\line(0,1){20}}
\put(0,80){\line(1,0){100}}

\put(0,60){\line(1,0){20}}
\put(0,20){\line(1,0){20}}
\put(9,35){$\vdots$}
\put(7,11){$p_s$}
\put(7,71){$p_2$}
\put(7,90){$t_1$}
\put(27,90){$t_2$}
\put(53,90){$\ldots$}
\put(87,90){$t_i$}
\put(-25,48){$C = $}
\put(100,80){\line(0,1){20}}
\put(0,100){\line(1,0){100}}
\end{picture}
\end{center} 
with $1 \leq t_1 < t_2 < \cdots < t_i \leq r$ and $t_1 \leq p_2 \leq 
\cdots \leq p_s \leq r$. To make this resolution graded, assign
a degree to each tableau basis element as follows:
\[\deg C :=  d_{t_1} + d_{t_2} + \cdots + d_{t_i} + 
d_{p_2} + \cdots+  d_{p_s}\] 
where $d_{i_j}$ is the $i_j$th element of $(d_1,\ldots,d_r)$.
\end{remark}

If $J$ is a homogeneous ideal of $R$, then the {\it Hilbert
function} of $R/J$ is $H_{R/J}(i) := \dim_k (R/J)_i$.
As a corollary of Theorem \ref{maintheorem} we
obtain a formula for the Hilbert function of
$R/I^s$.  We define
\[
\LL_{r,s} := \left\{ (a_1,\ldots,a_r) \in \N^r ~\left|~
a_1 + \cdots + a_r \leq s-1 \right\}\right..
\]
\begin{corollary} \label{hilbertfunction}
With the hypotheses as in Theorem \ref{maintheorem},
\[H_{R/I^s}(i) = \sum_{(a_1,\ldots,a_r) \in \LL_{r,s}}
H_{R/I}(i-a_1d_1-\cdots-a_rd_r)\]
where $H_{R/I}$ is the Hilbert function of $R/I$.
\end{corollary}

\begin{remark}
Corollary \ref{hilbertfunction} can also be proved directly
by using Lemma \ref{inductionstep} and an induction proof similar to the  proof
of Theorem \ref{maintheorem}.  The fact that  $H_{R/I^s}$ depends only upon
$H_{R/I}$ was first observed in \cite{CV}.
\end{remark}


\section{Application: A conjecture of Herzog, Huneke, and Srinivasan}

As an application of Theorem ~\ref{maintheorem}, we can verify
a special case of a conjecture of Herzog, 
Huneke, and Srinivasan \cite{HS}.  The following relation
between the multiplicity $e(R/I)$ of $R/I$ and the degrees of the syzygies
of $I$ is conjectured to hold:

\begin{conjecture} \label{HS}
Let $R/I$ be a Cohen-Macaulay $k$-algebra
with resolution of the form
\[0 \longrightarrow   \bigoplus_{j=1}^{b_r}
R(-d_{rj}) \longrightarrow \cdots \longrightarrow 
\bigoplus_{j=1}^{b_1}
R(-d_{1j}) \longrightarrow R \longrightarrow R/I
\longrightarrow 0.\]
Set $m_i = \min\{ d_{ij} ~|~ j = 1,\ldots,b_i\}$ and $M_i
= \max\{d_{ij} ~|~ j = 1,\ldots,b_i\}$.  Then 
\[\frac{\prod_{i=1}^r m_i}{r!} \leq e(R/I) \leq 
\frac{\prod_{i=1}^r M_i}{r!}.\]
\end{conjecture} 

By using Theorem \ref{maintheorem}, powers of complete
intersections can be added to the list of ideals in \cite{G,HS2} that satisfy
the conjecture.  The case $s=1$ is also found in \cite{HS}.

\begin{proposition}
Let $I$ be a complete intersection of $R$ of type $(d_1,\ldots,d_r)$,
and let $s$ be any positive integer.
Then Conjecture
\ref{HS} is true for $R/I^s$.
\end{proposition}

\begin{proof}
We can assume $d_1 \leq \cdots \leq d_r$.  
The resolution of $I^s$ given in Theorem \ref{maintheorem} implies
$m_i = sd_1 + d_2 + \cdots + d_i$ and  $M_i = d_{r-i+1}
+ \cdots + d_{r-1} + sd_r$
for $i = 1,\ldots,r$.  Because $e(R/I^s) = \binom{s+r-1}{r}d_1\cdots d_r$
and 
\begin{eqnarray*}
\frac{\prod^{r}_{i=1} m_i}{r!}\! & \!\leq &\!
\prod^{r}_{i=1} \frac{s+(i-1)}{i}d_i \! =\! e(R/I^s) =
 \prod^{r}_{i=1} \frac{s+(i-1)}{i}d_{r-i+1}
\leq \frac{\prod^{r}_{i=1} M_i}{r!},
\end{eqnarray*}
$R/I^s$ satisfies the conjecture.
\end{proof}


\section{Application:  Fat points with support on a complete intersection}

We apply Theorem \ref{maintheorem}  to the study of invariants
associated to sets of fat points in $\pr^n$ whose support is 
a complete intersection of points.  For this section and the next, we shall
assume that the field $k$ is also algebraically closed.

Let $P_1,\ldots,P_s$ be $s$ points in $\pr^n:=\pr^n_k$.  If $m_1,\ldots,m_s$
are $s$ positive integers, then let $Z = \{(P_1,m_1),\ldots,(P_s,m_s)\}$ 
denote the
subscheme of $\pr^n$ defined by
\[\Iz = \wp^{m_1}_1 \cap \wp^{m_2}_2 \cap \cdots \cap \wp^{m_s}_s\]
where $\wp_i$ is the defining ideal in $R$ of $P_i$.
We call $Z$ a {\it fat point scheme}, or a {\it set of
fat points.}  If $m_1 = \cdots = m_s = m$,
then we refer to $Z$ as a {\it homogeneous scheme of fat points},
otherwise $Z$ is {\it non-homogeneous}.  The set Supp$(Z) = \{P_1,\ldots,
P_s\}$ is the support of $Z$.

If $Z$ is a homogeneous fat point scheme of multiplicity $m$,
then $I_Z = I_X^{(m)}$, the $m$th symbolic power of $I_X$ where
$X = \operatorname{Supp}(Z)$.  For an arbitrary support, $I_X^m 
\subset I_X^{(m)} = I_Z$ since $I_X^m$ may not be saturated.  

If $I = (F_1,\ldots,F_n)$ is a complete intersection of type $(d_1,
\ldots,d_n)$, and if $I = \sqrt{I}$, then $I$ is the defining
ideal of $\prod_{i=1}^n d_i$ reduced points in $\pr^n$.  We denote this
set by $X = CI(d_1,\ldots,d_n)$ and call $X$ a complete intersection.
For each $m\in \N^+$ we write $Z = \{CI(d_1,\ldots,d_n);m\}$ to
denote the homogeneous fat point scheme
of multiplicity $m$ whose
support is $CI(d_1,\ldots,d_n)$. 
When $Z =\{CI(d_1,\ldots,d_n);m\}$, then we have
the equality 
$I_Z = I_X^{(m)} = I_X^m$ because of the following lemma:

\begin{lemma}[{\cite[Lemma 5, Appendix 6]{ZS}}]
If $I$ is a complete intersection, then $I^m = I^{(m)}$ for all
positive integers $m$.
\end{lemma}

Since the defining ideal of $\{CI(d_1,\ldots,d_n);m\}$ is a power
of a complete intersection, we therefore have

\begin{proposition}\label{hfci}
Suppose $Z = \{CI(d_1,\ldots,d_n);m\} \subseteq \pr^n$ with defining 
ideal $I_Z$. 
Then the graded minimal free resolution of $I_Z$ is given by 
Theorem \ref{maintheorem}.  The Hilbert function
$H_{R/I_Z}$
is given by Corollary \ref{hilbertfunction}.
\end{proposition}

If $Z$ is 
a fat point scheme with Supp$(Z) = CI(d_1,\ldots,d_n)$, but
$Z$ is not homogeneous, then Proposition \ref{hfci} can be used
to obtain partial information about the invariants associated to $Z$.
Recall that  $\alpha(Z) := \min\{i ~|~ (I_Z)_i \neq 0\}$ and 
\[
\ri(Z) := \min\left\{i ~\left|~ H_{R/I_Z}(i) = \deg(Z):
= \sum_{i=1}^s \binom{n+m_i-1}{n}\right\}\right.. 
\]
The invariant $\ri(Z)$ is 
the regularity index of $Z$, while $\alpha(Z)$ is the smallest
degree of a form contained in $I_Z$.

\begin{proposition}
Let $Z = \{(P_1,m_1),\ldots,(P_s,m_s)\} \subseteq \pr^n$ be a
non-homoge\-ne\-ous fat point scheme with Supp$(Z) =X = CI(d_1,\ldots,d_n)$.
Set $M = \max\{m_i\}_{i=1}^s$ and $m= \min\{m_i\}_{i=1}^s$.
Then
\[\sum_{\LL_{n,m}} H_{R/I_X}(i-a_1d_1-\cdots-a_nd_n) \leq
H_{R/I_Z}(i) \leq \sum_{\LL_{n,M}}  H_{R/I_X}(i-a_1d_1-\cdots-a_nd_n). \]
In particular
\begin{enumerate}
\item[$(i)$] $md_1 \leq \alpha(Z) \leq Md_1$.
\item[$(ii)$]$d_1+\cdots+d_{n-1} + md_n -n \leq \ri(Z) \leq
d_1 + \cdots + d_{n-1} + Md_n -n$.
\end{enumerate}
\end{proposition}

\begin{proof}
If $I_X$ is the defining ideal of the support $CI(d_1,\ldots,d_n)$,
then we have
\[ I_X^M\subseteq I_Z = \wp_1^{m_1} \cap \cdots \cap \wp_s^{m_s}
\subseteq I_X^m.\]
Since $I_X^M$ and $I_X^m$ 
define homogeneous fat point schemes on a complete intersection,
the conclusions now follow from Proposition \ref{hfci}.
\end{proof}

\section{Application:  Fat points with support on a complete 
intersection minus 
a point}

We discuss the invariants of fat point schemes whose support is a complete
intersection minus a point.  
To provide some motivation, we recall
that a set of points $X$
has the Cayley-Bacharach property (CBP)
if for every $P \in X$, $Y = X \backslash \{P\}$ always has
the same Hilbert function. It is well known (see \cite{GKR}) 
that a complete intersection
of points satisfies the CBP because of the following result:

\begin{theorem}\label{cbp}
Let $X = CI(d_1,\ldots,d_n) \subseteq \pr^n$ be a complete intersection.  
Let $P \in X$ be any point, and set $Y = X \backslash \{P\}$.  Then
\[H_Y(i) = \min\{H_X(i),~|X|-1\} ~~\mbox{for all $i$}.\]
\end{theorem}
Moreover, because $H_X$ depends only upon the type
$(d_1,\ldots,d_n)$, it follows that
the Hilbert function of $Y = X \backslash \{P\}$ also
depends upon the type.  
 
Since the Hilbert function of $Z = \{CI(d_1,\ldots,d_n);m\}$ depends
upon the type and $m$, it is natural to wonder if a 
``Cayley-Bacharach like'' result holds for $Z$, that is, 
if $(P,m)$ is any fat point of $Z$, does the Hilbert function
of $Y = Z \backslash\{(P,m)\}$  depend only upon the type and $m$?
As the next two examples show, the answer is no since the 
construction of the underlying complete intersection
must be taken into account. This suggests
that it may be difficult to find general results to describe
the invariants in this case.

\begin{example}  \label{example1}
Let $X = CI(3,4)$ be a complete intersection in $\pr^2$
with defining ideal $I_X = (F,G)$.
Set $Z = \{CI(3,4);3\}$.
Take any point $P \in X = \operatorname{Supp}(Z)$, and set
$Y = Z \backslash\{(P,3)\}$.
If $F = L_1L_2L_3$ and $G = L'_1L'_2L'_3L'_4$ are the product of linear forms,
then the Hilbert function of 
$Y$ is:
\begin{center}
\begin{tabular}{cccccccccccccccc}
$H_{Y}(i)$ & :  &1 &3 &6 &10 &15 &21 &28 &36 &45 &54 &62 &65
&66 &
$\rightarrow$ \\
\end{tabular}
\end{center}
On the other hand, if $F$ and $G$ are both irreducible, then
$H_Y$ is given by
\begin{center}
\begin{tabular}{cccccccccccccccc}
$H_{Y}(i)$& :  &1 &3 &6 &10 &15 &21 &28 &36 &45 &54 &62 &66 &
$\rightarrow$ &
\end{tabular}
\end{center}
The two Hilbert functions do not agree when $i=11$. 
\end{example}

\begin{example}  \label{example2}
This example shows that, in some cases, the Hilbert function depends
upon what point is removed from the underlying complete
intersection.
Let $G_1, G_2$ be two irreducible
forms of $R = k[x,y,z]$ with $\deg G_1 = 2$ and $\deg G_2 =
3$. Let $L_1$ and $L_2$ be two linear forms that do not pass
through the points in the complete intersection defined by
$(G_1,G_2)$. Set $F_1 = L_1G_1$ and $F_2 = L_2G_2$, and let $X =
CI(3,4)$ denote the complete intersection of $\pr^2$ defined by $I_X =
(F_1,F_2)$. Let $P_1$ denote the point of $X$  with defining
ideal $I_{P_1} = (L_1,L_2)$, and  let $P_2 \in X$ be a point that
lies in the complete intersection defined by $(G_1,G_2)$.  Let
$Y_1 = \{X;3\} \backslash \{(P_1,3)\}$ and
$Y_2 = \{X;3\} \backslash \{(P_2,3)\}$. 
We then have
\begin{center}
\begin{tabular}{ccccccccccccccccc}
$H_{Y_1}(i)$ & :  &1 &3 &6 &10 &15 &21 &28 &36 &45 &54 &62 &65
&66 &
$\rightarrow$ \\
$H_{Y_2}(i)$& :  &1 &3 &6 &10 &15 &21 &28 &36 &45 &54 &62 &66 &
 66 & $\rightarrow$ &
\end{tabular}
\end{center}
which fail to agree when $i=11$.
\end{example}

\begin{remark} In her Ph.D thesis \cite{Gu} the first author 
introduced a generalization of the Cayley-Bacharach property for fat points.
Specifically, a homogeneous set of fat points $Z \subseteq \pr^n$ 
is said to satisfy the generalized
Cayley-Bacharach property if for all fat points ${(P,m)} \in Z$,
the set of fat points $Y = Z \backslash \{(P,m)\}$  has the same Hilbert
function.  
As Example \ref{example2} shows, 
even if the support is a complete intersection, 
this generalized property need not hold.
\end{remark}

As noted, these examples imply that extra hypotheses are needed on
the support in order to study fat points whose support is a complete
intersection minus a point.  We therefore end this paper by 
considering $\alpha(Y)$
and $\ri(Y)$ of $Y = \{X;m\} \backslash \{(P,m)\}$ when the
support $X$ ``splits'' into smaller complete intersections,
say $X = X_1 \cup X_2$, with the property that the point $P$ being
removed from $X$ belongs to $X_2$ and $X_2$ is contained in a 
hyperplane, i.e, $X_2 \subseteq H \cong\pr^{n-1} \subseteq \pr^n$.  
This class of support allows us to make further 
use of Theorem \ref{maintheorem}.
We introduce some relevant terminology.

\begin{definition}  The complete intersection $X = CI(d_1,\ldots,d_n)
\subseteq \pr^n$ {\it splits} if $X$ can be 
written as the the disjoint union of two complete intersections:
\[X = CI(d_1,\ldots,d_{i-1},d_i-1,d_{i+1},\ldots,d_n)
\cup CI(d_1,\ldots,d_{i-1},1,d_{i+1},\ldots,d_n)\]
for some $1 \leq i \leq n$ with $d_i \geq 2$. 
\end{definition}

\begin{theorem}  \label{alphabound2}
Let $X = CI(d_1,\ldots,d_n) \subseteq \pr^n$
with $d_1 \leq \cdots \leq d_n$,
and suppose that $X$ splits as $X = CI(d_1,\ldots,d_n-1) \cup
CI(d_1,\ldots,d_{n-1},1)$.
Let $Z = \{X;m\}$ and $Y = \{X;m\} \backslash \{(P,m)\}$
for any $P \in CI(d_1,\ldots,d_{n-1},1) \subseteq X$. 
\begin{enumerate}
\item[$(i)$] If $d_1 = d_n$, then $m(d_1-1) \leq \alpha(Y) \leq
\alpha(Z) = md_1$.
\item[$(ii)$]If $d_1 < d_n$, then  $\alpha(Y) = \alpha(Z) = md_1$.
\end{enumerate}
\end{theorem}

\begin{proof}
The complete intersection $CI(d_1,\ldots,d_{n-1},d_n-1)$
is a subset of Supp$(Y)$ $= X \backslash \{P\}$.  Hence
$Y' = \{CI(d_1,\ldots,d_{n-1},d_n -1);m\}$ is a subscheme
of $Y$.  This implies $I_Y \subseteq I_{Y'}$, and thus
 $\alpha({Y'}) \leq \alpha(Y)$.  But
since the support of $Y'$ is a complete intersection, by 
Theorem \ref{hfci} we have
$\alpha({Y'}) = m\cdot\min\{d_1,\ldots,d_n-1\}$.
Because $d_1 \leq \cdots \leq d_n$,
if $d_1 = d_n$, then $\min\{d_1,\ldots,d_n-1\} = d_n-1$,
thus proving $(i)$.  In the case of $(ii)$,
$\min\{d_1,\ldots,d_n-1\} = d_1$.
\end{proof}

If $Z \subseteq \pr^n$ is any fat point scheme, and if $Y \subseteq Z$, 
then Lemma 1.1 of \cite{TV}
implies $\ri(Y) \leq \ri(Z)$.  Thus, if 
$Y \subseteq Z = \{CI(d_1,\ldots,d_n);m\}$, we have by Theorem
\ref{hfci}
\[\ri(Y) \leq \ri(Z) = d_1 + \cdots + d_{n-1} + md_n -n.\]

\begin{proposition} \label{lowerbound-splittable}
Let $X = CI(d_1,\ldots,d_n) \subseteq \pr^n$
with $2 \leq d_1 \leq \cdots \leq d_n$ and suppose $X$ splits as
$X = X_1 \cup X_2 = CI(d_1-1,d_2,\ldots,d_n) \cup
CI(1,d_2,\ldots,d_n)$.
If $Z = \{X;m\}$, and $Y = Z \backslash \{(P,m)\}$ with
$P \in X_2$, then
\[\ri(Y) \geq d_1 + \cdots + d_{n-1} + md_n - (n+1).\]
\end{proposition}

\begin{proof}  Note that $Z = Z_1 \cup Z_2$ where $Z_1 = \{X_1;m\}$ and
$Z_2 = \{X_2;m\}$.
Hence $Y = Z_1 \cup Y_2$ where $Y_2 =\{X_2;m\} \backslash \{(P,m)\}$.
From the short exact sequence
\[0 \rightarrow R/I_Y \rightarrow R/I_{Z_1} \oplus R/I_{Y_2}
\rightarrow R/(I_{Z_1} + I_{Y_2}) \rightarrow 0\]
it follows that
$H_{Y}(i) \leq H_{Z_1}(i) + H_{Y_2}(i)$ for all $i \in \N$.
Set $r = d_1 + \cdots + d_{n-1} + md_n - (n+2)$.  Then
$H_{Z_1}(r) < \deg Z_1 = (d_1-1)d_2\cdots d_n \binom{m+n-1}{n}$
by Theorem \ref{hfci} since $Z_1$ is a homogeneous fat point
scheme on a complete intersection.
Hence
$H_Y(r) < \deg Z_1 + \deg Y_2
= \deg Y,$ and thus $\ri(Y) > r$.
\end{proof}
Under the hypotheses of the previous proposition we have
\[d -(n+1) \leq \ri(Y) \leq d - n\]
where $d = d_1 + \cdots + d_{n-1} +md_n.$ If
$X = CI(d_1,d_2) \subseteq \pr^2$, we can give an exact formula
when the support splits nicely.
If $X \subseteq \pr^2$ is any set of points, following \cite{TV}, we define
$b(X) := \min\{t ~|~ I_t
~\mbox{contains a regular sequence of length two}\}.$
A bound on the regularity index for any set of fat points in $\pr^2$ is
then given by 

\begin{theorem}[{\cite[Theorem 3.2]{TV}}]  \label{trungvalla}
Let $X = \{P_1,\ldots,P_s\}$
be a set of points in $\pr^2$ with $m_1 \geq \cdots \geq m_s$ a sequence
of positive integers.  For $i = 1,\ldots,m_1$, let
$Y_i = \{P_j ~|~ m_j \geq i\}.$
If $Z = \{(P_1,m_1),\ldots,(P_s,m_s)\}$, then
\[\ri(Z) \leq \ri(X) + \sum_{i=2}^{m_1} b(Y_i).\]
\end{theorem}
If $X = CI(d_1,d_2)$, then $b(X) = d_2$.
Furthermore, since $(\Ix)_t = (\Iy)_t$ for $t \leq d_1 + d_2 -3$
by Theorem \ref{cbp},
we also have $b(Y) = d_2$.

\begin{corollary}\label{stessobound}
Let $X = CI(d_1,d_2) \subseteq \pr^2$
with $2 \leq d_1 \leq d_2$ and suppose $X$ splits as
$X = X_1 \cup X_2 = CI(d_1-1,d_2) \cup
CI(1,d_2)$.
If $Z = \{X;m\}$, and $Y = Z \backslash \{(P,m)\}$ with
$P \in X_2$, then
\[\ri(Y) = d_1 + md_2 - 3.\]
\end{corollary}

\begin{proof}
For $i = 1,\ldots,m$, $Y_i = X\backslash \{P\}$ and therefore,
$b(Y_i) = d_2$.  The formula of Theorem \ref{trungvalla}
and the fact
that $\ri(X\backslash \{P\}) = d_1 + d_2 -3$ implies $\ri(Y) \leq
d_1 + md_2 - 3.$  The reverse inequality is Proposition \ref{lowerbound-splittable}.
\end{proof}

\begin{remark}
If $I_X = (F_1,\ldots,F_n)$ is the defining ideal of the
complete intersection $X = CI(d_1,\ldots,d_n)
\subseteq \pr^n$, and if each $F_i$ is the product of $d_i$
distinct linear forms, then the points of $X$ form a $n$-dimensional ``box''.
For each $i = 1,\ldots,n$, $X$ splits as 
\[X = CI(d_1,\ldots,d_i-1,\ldots,d_n) \cup CI(d_1,\ldots,1,\ldots,d_n).\]
Hence, the previous results about fat points on complete intersections
which split apply to these configurations.
\end{remark}



\end{document}